\documentclass[12pt]{article}

\usepackage{graphicx,tikz,color,url}
\usepackage{amsmath,amssymb,latexsym}
\usepackage{MnSymbol}

\newtheorem{theorem}{Theorem}[section]
\newtheorem{definition}[theorem]{Definition}
\newtheorem{lemma}[theorem]{Lemma}
\newtheorem{proposition}[theorem]{Proposition}
\newtheorem{observation}[theorem]{Observation}
\newtheorem{conjecture}[theorem]{Conjecture}
\newtheorem{corollary}[theorem]{Corollary}

\newtheorem{claim}{Claim}

\newcommand{\proof}{\noindent{\bf Proof.\ }}
\newcommand{\qed}{\hfill $\square$ \bigskip}

\newcommand{\cF}{{\cal F}}
\newcommand{\cU}{{\cal U}}
\newcommand{\cO}{{\cal O}}
\newcommand{\cQ}{{\cal Q}}
\newcommand{\cV}{{\cal V}}

\newcommand{\smallqed}{{\tiny ($\Box$)}}

\newcommand{\gtg}{\gamma_{tg}}

\begin{document}

\title{Perfect graphs for domination games}

\author{Csilla Bujt\'as$^{a}$\thanks{Email: \texttt{csilla.bujtas@fmf.uni-lj.si}} 
\and Vesna Ir\v si\v c$^{a,b}$\thanks{Email: \texttt{vesna.irsic@fmf.uni-lj.si}}
\and Sandi Klav\v zar $^{a,b,c}$\thanks{Email: \texttt{sandi.klavzar@fmf.uni-lj.si}}
}
\maketitle

\begin{center}
$^a$ Faculty of Mathematics and Physics, University of Ljubljana, Slovenia\\
\medskip

$^b$ Institute of Mathematics, Physics and Mechanics, Ljubljana, Slovenia\\
\medskip

$^c$ Faculty of Natural Sciences and Mathematics, University of Maribor, Slovenia\\
\medskip
\end{center}

\begin{abstract}
Let $\gamma_g(G)$ and $\gamma_{tg}(G)$ be the game domination number and the total game domination number of a graph $G$, respectively. Then $G$ is \emph{$\gamma_g$-perfect} (resp.\ \emph{$\gamma_{tg}$-perfect}), if every induced subgraph $F$ of $G$ satisfies $\gamma_g(F)=\gamma(F)$ (resp.\ $\gamma_{tg}(F)=\gamma_t(F)$). A recursive characterization of $\gamma_g$-perfect graphs is derived. The characterization yields a polynomial recognition algorithm for $\gamma_g$-perfect graphs. It is proved that every minimally $\gamma_g$-imperfect graph has domination number $2$. All minimally $\gamma_g$-imperfect triangle-free graphs are determined. It is also proved that $\gamma_{tg}$-perfect graphs are precisely $\overline{2P_3}$-free cographs. 
\end{abstract}


\noindent
{\bf AMS Subj.\ Class.\ (2010)}: 05C57, 05C69, 68Q25

%%%%%%%%%%%%%%%%%%%%%%%%%%%%%%%%%%%%%%%%%%%%%%%%%%%%%
%%%%%%%%%%%%%%%%%%%%%%%%%%%%%%%%%%%%%%%%%%%%%%%%%%%%%
\section{Introduction}
\label{sec:intro}
%%%%%%%%%%%%%%%%%%%%%%%%%%%%%%%%%%%%%%%%%%%%%%%%%%%%%
%%%%%%%%%%%%%%%%%%%%%%%%%%%%%%%%%%%%%%%%%%%%%%%%%%%%%

The domination game on a graph $G$ is played by Dominator and Staller. If Dominator (resp.\ Staller) starts the game, we speak of the D-game (resp.\ S-game). During the game, the players alternatively select vertices that are not dominated by the set of previously selected vertices. The game ends when no such vertex is available. Dominator's goal is to finish the game as soon as possible, while Staller wishes to play the game as long as possible. The unique number of moves played in the D-game (resp.\ S-game) when both players play optimally is the game domination number $\gamma_g(G)$ (resp.\ Staller-start game domination number $\gamma_g'(G)$) of $G$. The total domination game is defined analogously, the only difference being that when a new vertex is selected, it must totally dominate at least one vertex not yet totally dominated by the previously selected vertices; the corresponding game total domination numbers are denoted with $\gamma_{tg}(G)$ and $\gamma_{tg}'(G)$.   

The seminal paper~\cite{bresar-2010} on the domination game together with its follow-up~\cite{kinnersley-2013} had a great impact, leading to several dozens of papers. Instead of listing (too) long list of references, we just point to~\cite{dorbec-2015, james-2019, klavzar-2019, xu-2018} and references therein. Similarly, the seminal papers on the total domination game~\cite{henning-2015, henning-2017} led to its extensive investigation, cf.~\cite{bujtas-2018, bujtas-2016, henning-2018}. 

Perfect graphs lie in the very core of graph theory, papers~\cite{chudnovsky-2006, lovasz-1972} being highlights of the theory. Now, just as $\chi(G)\ge \omega(G)$ holds trivially for every graph $G$, we infer from definitions that $\gamma_g(G)\ge \gamma(G)$. Hence we say that $G$ is a \emph{$\gamma_g$-minimal graph} if the equality $\gamma_g(G)=\gamma(G)$ holds. In this paper, we will call these graphs \emph{$\gamma_g$-graphs} for short. The $\gamma_g$-minimal trees were characterized in~\cite{nadjafi-2016}, but a characterization of $\gamma_g$-minimal graphs is widely open. In this paper we study the hereditary version of this property via the following concept.  

\begin{definition}
A graph $G$ is \emph{$\gamma_g$-perfect}, if every induced subgraph $F$ of $G$ satisfies $\gamma_g(F)=\gamma(F)$.
\end{definition}
Note that $G$ and/or $F$ may be disconnected in the above definition.

Since the inequalities $\gamma_g'(G)\ge \gamma(G)$, $\gamma_{tg}(G)\ge \gamma_t(G)$, and $\gamma_{tg}'(G)\ge \gamma_t(G)$  also hold for every graph $G$ (where, if the total domination is involved, $G$ must of course be isolate-free), we introduce the analogous terminology for the Staller-start domination game and for the total domination games.
\begin{itemize}
 \item $G$ is a \emph{$\gamma_g'$-minimal graph} if $\gamma(G)=\gamma_g'(G)$ holds, and is \emph{$\gamma_g'$-perfect} if all of its induced subgraphs are $\gamma_g'$-minimal graphs. 
 \item $G$ is a \emph{$\gtg$-minimal graph} if $\gamma_t(G)=\gamma_{tg}(G)$ holds, and is \emph{$\gamma_{tg}$-perfect} if all of its isolate-free induced subgraphs are $\gtg$-minimal graphs. 
 \item $G$ is a \emph{$\gtg'$-minimal graph} if $\gamma_t(G)=\gamma_{tg}'(G)$ holds, and is \emph{$\gamma_{tg}'$-perfect} if all of its isolate-free induced subgraphs are $\gtg'$-minimal graphs.
\end{itemize}

In the literature, several similar problems were studied, that is, the equality between two covering and/or domination-type invariants is required to hold not only for a graph $G$ but also for all induced subgraphs. See~\cite{SM, ZZ} for earlier approaches and~\cite{ADR, ABBT, BB-2018, CP, HJR} for recently published results. 

In this paper we characterize $\gamma_g$-perfect, $\gamma_g'$-perfect, $\gamma_{tg}$-perfect, and $\gamma_{tg}'$-perfect graphs. The main result is a characterization of the first class, it is proved in Section~\ref{sec:gamma-g}. The characterization describes a recursive structure of $\gamma_g$-perfect graphs that in particular yields a polynomial recognition algorithm for $\gamma_g$-perfect graphs. We also introduce minimally $\gamma_g$-imperfect graphs and prove that they have domination number $2$. In Section~\ref{sec:appl} we discuss recognition complexity of $\gamma_g$-perfect graphs, and present some results on triangle-free graphs. Characterizations of $\gamma_g'$-perfect, $\gamma_{tg}$-perfect, and $\gamma_{tg}'$-perfect turned out to be simpler, we state them in Section~\ref{sec:further}. In particular, $\gamma_{tg}'$-perfect graphs are precisely cographs, and $\gamma_{tg}$-perfect graphs are precisely $\overline{2P_3}$-free cographs. 

%%%%%%%%%%%%%%%%%%%%%%%%%%%%%%%%%%%%%%%%%%%%%%%%%%%%%
%%%%%%%%%%%%%%%%%%%%%%%%%%%%%%%%%%%%%%%%%%%%%%%%%%%%%
\section{Preliminaries}
\label{sec:preliminaries}
%%%%%%%%%%%%%%%%%%%%%%%%%%%%%%%%%%%%%%%%%%%%%%%%%%%%%
%%%%%%%%%%%%%%%%%%%%%%%%%%%%%%%%%%%%%%%%%%%%%%%%%%%%%

If $v$ is a vertex of a graph $G=(V(G), E(G))$, then the \emph{open neighborhood} $N_G(v)$ is the set of neighbors of $v$, while the \emph{closed neighborhood} $N_G[v]$ is the open neighborhood supplemented with the vertex $v$ itself. Two vertices, $u$ and $v$, are \emph{(true) twins} in $G$, if $N_G[u]=N_G[v]$, and they are \emph{false twins} if $N_G(u)=N_G(v)$. The \emph{degree} of $v$ in $G$ is $d_G(v)=|N_G(v)|$. The \emph{closed neighborhood of a set $S$} of vertices is $N_G[S]=\bigcup_{v \in S}N_G[v]$. In this paper, the \emph{open neighborhood of $S$} will be meant as $N'_G(S)= N_G[S]\setminus S$.

A set $S \subseteq V(G)$ is a \emph{dominating set} of $G$ if $N_G[S] = V(G)$. The minimum cardinality of a dominating set is the \emph{domination number} $\gamma(G)$ of $G$. A dominating set $S$ is a \emph{total dominating set} if every vertex from $S$ has a neighbor in $S$. The smallest cardinality of a total dominating set is the \emph{total domination number} $\gamma_{t}(G)$ of $G$, see the book~\cite{MHAYbookTD}. 

The distance $d_G(u,v)$ between vertices $u$ and $v$ of a connected graph $G$ is the minimum number of edges on a $u,v$-path. If $H_1$ and $H_2$ are subgraphs of a connected graph $G$, then the distance $d_G(H_1,H_2)$ between $H_1$ and $H_2$ is the minimum of the distances $d_G(v_1, v_2)$, where $v_1\in V(H_1)$ and $v_2\in V(H_2)$. 

We will say that a graph $G$ is \emph{minimally $\gamma_g$-imperfect}, if each of its proper induced subgraphs is $\gamma_g$-perfect but $\gamma(G) < \gamma_g(G)$. Since perfectness is a hereditary property, a graph is not $\gamma_g$-perfect if and only if it has an induced subgraph which is minimally $\gamma_g$-imperfect. This ensures that there exists a forbidden subgraph characterization for $\gamma_g$-perfect graphs.    \emph{Minimally $\gamma_g'$-imperfect}, \emph{minimally $\gamma_{tg}$-imperfect} and \emph{minimally $\gamma_{tg}'$-imperfect} graphs are defined analogously to minimally $\gamma_g$-imperfect graphs.

We say that a graph $G$ is \emph{$2$-$\gamma_g$-perfect}, if every induced subgraph $F$ of $G$ with $\gamma(F)=2$ is a $\gamma_g$-graph. By definition, every induced subgraph of a $2$-$\gamma_g$-perfect graph is $2$-$\gamma_g$-perfect, as well. This concept will be useful when proving our characterization theorem for $\gamma_g$-perfect graphs, but at the end we show that this property is equivalent to the $\gamma_g$-perfectness.

{\em Cographs} are, by definition, the graphs that contain no induced path $P_4$. These graphs admit different characterizations, see~\cite{corneil-1981};  the one to be applied here asserts that cographs are precisely the graphs that can be obtained from $K_1$ by means of the disjoint union and join of graphs. 

%%%%%%%%%%%%%%%%%%%%%%%%%%%%%%%%%%%%%%%%%%%%%%%%%%%%%
%%%%%%%%%%%%%%%%%%%%%%%%%%%%%%%%%%%%%%%%%%%%%%%%%%%%%
\section{Characterization of $\gamma_g$-perfect graphs}
\label{sec:gamma-g}
%%%%%%%%%%%%%%%%%%%%%%%%%%%%%%%%%%%%%%%%%%%%%%%%%%%%%
%%%%%%%%%%%%%%%%%%%%%%%%%%%%%%%%%%%%%%%%%%%%%%%%%%%%%

In this section our goal is to characterize $\gamma_g$-perfect graphs. 
The theorem, that is formulated and proved in Subsection~\ref{subsec:char}, states an equivalence with a recursively defined graph class. The two operators used in the recursive definition are introduced in Subsection~\ref{subsec:pre}. As a consequence of the characterization theorem, we prove in Subsection~\ref{subsec:recognition} that $\gamma_g$-perfect graphs can be recognized in polynomial time.

Along the proof of the main theorem, we first consider $2$-$\gamma_g$-perfect graphs and prove that they can be built from an isolated vertex by using the two specified operators. We also show that these operators applied to $2$-$\gamma_g$-perfect graphs always result in $\gamma_g$-graphs. Then, using further statements on the structure of $2$-$\gamma_g$-perfect graphs, we can prove the equivalence between $\gamma_g$-perfectness, $2$-$\gamma_g$-perfectness, and the property of recursive constructability.

%%%%%%%%%%%%%%%%%%%%%%%%%%%%%%%%%%%%%%%%%%%%%%%%%%%%%
\subsection{Preliminary observations}
\label{subsec:pre}
%%%%%%%%%%%%%%%%%%%%%%%%%%%%%%%%%%%%%%%%%%%%%%%%%%%%%

A \emph{homogeneous clique} $Q$ in a graph $G$ is a clique in which every two vertices are true twins, that is, $N_G[u] = N_G[v]$ holds for every $u,v \in V(Q)$. In other words, a clique $Q$ is homogeneous if there is a join between $Q$ and $N_G'(V(Q))$. From now on, we will use the same notation $Q$ when referring to the vertex set of the clique $Q$. A \emph{maximal homogeneous clique} (shortly, MHC) is an inclusion-wise maximal homogeneous clique.  By definition, two maximal homogeneous cliques are always vertex disjoint. A \emph{perfect set of cliques} (shortly, PSC) in a graph $G$ is a (possibly empty) set $\cQ$ of homogeneous cliques such that $d_G(Q, Q') = 3$ and there is a join between  $N'_G(Q)$ and $N'_G(Q')$ for every $Q, Q'\in \cQ$. Note that the empty set $\cQ= \emptyset$ is a PSC in $G$, and if $Q$ is a homogeneous clique in $G$, then the one-element set $\cQ=\{Q\}$ is also a PSC in $G$. If $\cQ$ is a nonempty PSC, we will usually use the notation $\cQ=\{Q_1,\dots, Q_k\}$ and $V(\cQ)=\bigcup_{i=1}^k Q_i$. For every PSC $\cQ$ which consists of at least two cliques, every $v\in N'_G(Q_i)$ has a neighbor which is not adjacent to the vertices of $Q_i$ and hence, the following statement holds by definitions.

\begin{observation}
	If $\cQ$ is a perfect set of cliques in a graph $G$ and $|\cQ|\ge 2$, then every $Q_i \in \cQ$ is a maximal homogeneous clique in $G$. 
\end{observation}

Given a graph $G$, its \emph{MHC-contraction} is the graph $\widehat{G}$ obtained from $G$ by contracting every maximal homogeneous clique into one vertex. Equivalently, a graph isomorphic to $\widehat{G}$ is obtained from $G$ by sequentially deleting one of two true twins while such a pair exists. For a vertex $v \in V(G)$, the corresponding vertex in $\widehat{G}$ will be denoted by $\widehat{v}$; that is, for any two vertices $u,v \in V(G)$, the vertices $\widehat{u}$ and $\widehat{v}$ are identical in $\widehat{G}$, if and only if $u$ and $v$ are true twins in $G$. For any two non-twin vertices $u,v\in V(G)$, by definition, we have $d_G(u,v)=d_{\widehat{G}}(\widehat{u},\widehat{v})$. Observe that if $Q$ is a homogeneous clique in $G$ and a vertex $v \in Q$ is dominated by a set $D \subseteq V(G)$, then every vertex $u\in Q$ is dominated by $D$. Similarly, if $v,u \in Q$ and $v\in D$, then $(D \cup \{u\}) \setminus \{v\}$ and $D$ dominate the same set of vertices. It follows that $|D\cap Q| \le 1$ holds for every minimum dominating set $D$ and homogeneous clique $Q$. We also infer the following facts.

\begin{observation}
	\label{observation}
	Every graph $G$ and its MHC-contraction $\widehat{G}$ satisfy the following statements:
	\begin{itemize}
		\item[$(i)$] $\gamma(G)=\gamma(\widehat{G})$;
		\item[$(ii)$] $\gamma_g(G)=\gamma_g(\widehat{G})$;
		\item[$(iii)$] $\gamma_g'(G)=\gamma_g'(\widehat{G})$;
		\item[$(iv)$] $G$ is a $\gamma_g$-perfect graph if and only if $\widehat{G}$ is  $\gamma_g$-perfect;
		\item[$(v)$] $G$ is a $\gamma_g'$-perfect graph if and only if $\widehat{G}$ is  $\gamma_g'$-perfect.
	\end{itemize}
\end{observation}

We will refer to the following graph operators:
\begin{itemize}
	\item Given a graph $G$ and a positive integer $s$, $G \bigcupdot K_s$ denotes the vertex disjoint union of $G$ and the complete graph $K_s$. 
   \item If $G$ is a graph, $v$ a new vertex and $\cQ$ a perfect set of cliques in $G$, then the graph $\cO(G, v, \cQ)$ is obtained from $G$ by adding the vertex $v$ and making it adjacent to all vertices in $V(G)\setminus V(\cQ)$.
\end{itemize}

%%%%%%%%%%%%%%%%%%%%%%%%%%%%%%%%%%%%
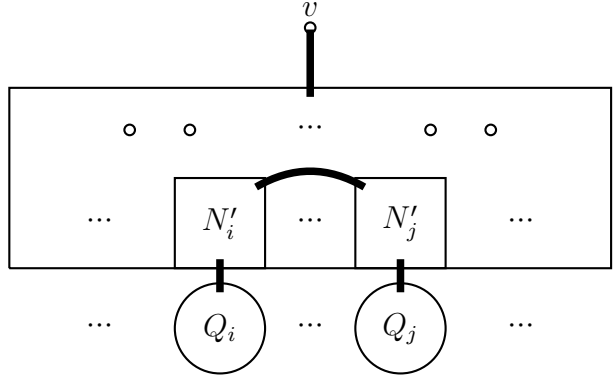
\begin{figure}[!ht]
	\begin{center}
		\begin{tikzpicture}[thick,scale=0.8]
		\draw (0,0) -- (10,0) -- (10,3) -- (0,3) -- (0,0);
		
		\draw (3.5, -1) circle (0.75cm);
		\node at (3.5, -1) {$Q_i$};
		
		\draw (6.5, -1) circle (0.75cm);
		\node at (6.5, -1) {$Q_j$};
		
		\node at (1.5, -1) {$\cdots$};
		\node at (5, -1) {$\cdots$};
		\node at (8.5, -1) {$\cdots$};
		
		\draw (2.75, 0) rectangle (4.25, 1.5);
		\node at (3.5, 0.75) {$N_i'$};
		
		\draw (5.75, 0) rectangle (7.25, 1.5);
		\node at (6.5,0.75) {$N'_j$};
		
		\node at (1.5, 0.75) {$\cdots$};
		\node at (5, 0.75) {$\cdots$};
		\node at (8.5, 0.75) {$\cdots$};
		
		\node[circle, draw, inner sep=0pt, minimum width=4pt] at (3,2.3) {};
		\node[circle, draw, inner sep=0pt, minimum width=4pt] at (2,2.3) {};
		\node[circle, draw, inner sep=0pt, minimum width=4pt] at (8,2.3) {};
		\node[circle, draw, inner sep=0pt, minimum width=4pt] at (7,2.3) {};
		\node at (5, 2.3) {$\cdots$};
		
		\node[circle, draw, inner sep=0pt, minimum width=4pt] at (5,4) {};

\draw (5,4.3) node {$v$};		

		\draw[line width=1mm] (3.5,-0.4) -- (3.5,0.15);
		\draw[line width=1mm] (6.5,-0.4) -- (6.5,0.15);
		\draw[line width=1mm] (5,3.97) -- (5,2.85);
		\draw[line width=1mm] (4.1, 1.35) to[out=30,in=150] (5.9, 1.35);
		
		\end{tikzpicture}
		\caption{An illustration of the operator $\cO(G, v, \cQ)$. The thick lines represent joins between the connected sets, $N_i'=N'_{G}(Q_i)$, and $N_j'=N'_{G}(Q_j)$.}
		\label{fig:transformation}
	\end{center}
\end{figure}

%%%%%%%%%%%%%%%%%%%%%%%%%%%%%%%%%%%%%
When referring to the operators, we will always assume that $s\in \mathbb{N}$ and $\cQ$ is a PSC in $G$. Note that, under this assumption,  $\gamma(G) \ge |\cQ|$ must be true as $d_G(Q_i, Q_j) > 2$ for every homogeneous cliques $Q_i$ and $ Q_j$ from $\cQ$ and therefore, the domination of the entire $V(\cQ)$ needs at least $|\cQ|$ different vertices in $G$. 

In Subsection~\ref{subsec:char} we will show that every $\gamma_g$-perfect graph can be built from an isolated vertex by using these two operators. In particular, we will prove that, if $G$ is $\gamma_g$-perfect, then both  $G \bigcupdot  K_s$ and $\cO(G, v, \cQ)$ are $\gamma_g$-perfect.
As a preliminary result, we show that under the stronger condition $\gamma(G) >|\cQ|$ the operator $\cO$ always gives a $\gamma_g$-graph even if $G$ is not $\gamma_g$-perfect.

\begin{proposition}
	\label{prop:opO}
 If $\cQ$ is a perfect set of cliques in $G$ and $\gamma(G) >|\cQ|$, then the graph $G'=\cO(G, v, \cQ)$ is a $\gamma_g$-graph.
  \end{proposition}
\proof If $\cQ=\emptyset$, then $v$ is a universal vertex  and $\gamma(G')=\gamma_g(G')=1$. Otherwise, let $\cQ=\{Q_1, \dots ,Q_k\}$. Choosing one vertex $x_i\in Q_i$ for every $i \in [k]$, we observe that $v,x_1,\dots x_k$ form a dominating set in $G'$ and therefore, $\gamma(G')\le k+1$. On the other hand, any dominating set $D$ which contains $v$, must contain a vertex $d_i$, which is different from $v$, to dominate $x_i$ for every $i\in [k]$. If $j\neq \ell$, then $d_G(x_j, x_\ell) = 3$ and hence, $d_j\neq d_\ell$. This proves $|D| \ge k+1$ if $v\in D$. If $v\notin D$, then $D$ is a dominating set also in $G$ and, by our condition, $|D|\ge \gamma(G) \ge k+1$. Therefore, $\gamma(G')=k+1$ holds.

In the domination game, let Dominator play $v$ as his first move. In the later turns, no matter how both players play, the domination of the homogeneous cliques $Q_1,\dots,Q_k$ needs exactly $k$ further vertices.  This strategy of Dominator shows that $\gamma_g(G')\le k+1=\gamma(G')$. Since $\gamma_g(G')\ge \gamma(G')$ is also true, we conclude that $\gamma_g(G')= \gamma(G')=k+1$. \qed

Extending a graph $G$ with a universal vertex, that is, constructing $\cO(G,v, \emptyset)$, always results in a $\gamma_g$-graph. Already this simple fact shows that any graph can be embedded into a $\gamma_g$-graph and consequently, the class of $\gamma_g$-graphs does not admit a forbidden subgraph characterization.  

In Subsection~\ref{subsec:char}, we will often use the following lemma that gives characterizations of   $\gamma_g$-graphs with small domination number. Each of these statements was either observed in \cite{KKS-2016} or can be obtained as a direct consequence of the earlier statements. 

\begin{proposition} 
	\label{lemma:1-2}
The following statements hold for every graph $G$.
	\begin{itemize}
		\item[$(i)$] $\gamma(G)=\gamma_g(G)=1$ holds if and only if $\gamma(G)=1$.  
		\item[$(ii)$] $\gamma(G)=\gamma_g(G)=2$ holds if and only if $G$ does not have a universal vertex but there exists a vertex $v\in V(G)$ such that $V(G)\setminus N\left[v \right]$ induces a homogeneous clique in $G$.
		\item[$(iii)$] If  $\Delta(G)= |V(G)|-2$, then $\gamma(G)=\gamma_g(G)=2$.
		\item[$(iv)$] If $G$ does not contain true twins, then  $\gamma(G)=\gamma_g(G)=2$ is true if and only if $\Delta(G)= |V(G)|-2$.
			\end{itemize}
\end{proposition}

%%%%%%%%%%%%%%%%%%%%%%%%%%%%%%%%%%%%%%%%%%%%%%%%%%%%%
\subsection{Minimally $\gamma_g$-imperfect graphs}
\label{subsec:imperfect}
%%%%%%%%%%%%%%%%%%%%%%%%%%%%%%%%%%%%%%%%%%%%%%%%%%%%%

In this subsection, we identify a collection of minimally $\gamma_g$-imperfect graphs. 

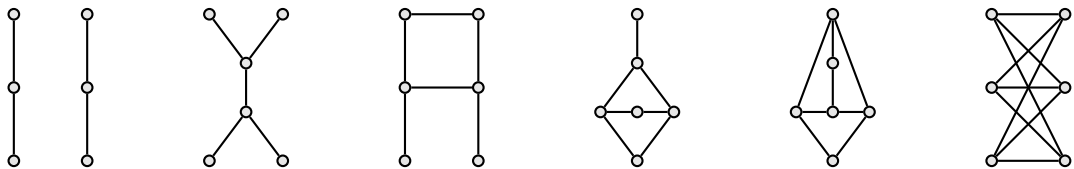
\begin{figure}[!ht]
	\begin{center}
		\begin{tikzpicture}[thick,scale=0.65]
		
		% Define style for nodes
		\tikzstyle{every node}=[circle, draw, fill=black!10,
		inner sep=0pt, minimum width=4pt]
		
		\begin{scope}
		\node (1) at (0, 0) {};
		\node (2) at (0, 1.5) {};
		\node (3) at (0, 3) {};
		\node (4) at (1.5, 0) {};
		\node (5) at (1.5, 1.5) {};
		\node (6) at (1.5, 3) {};
		
		\path (1) edge (2);
		\path (3) edge (2);
		\path (4) edge (5);
		\path (5) edge (6);
		\end{scope}
		
		\begin{scope}[xshift= 4cm]
		\node (1) at (0, 0) {};
		\node (2) at (1.5, 0) {};
		\node (3) at (0.75, 1) {};
		\node (4) at (0.75, 2) {};
		\node (5) at (0, 3) {};
		\node (6) at (1.5, 3) {};
		
		\path (1) edge (3);
		\path (3) edge (2);
		\path (4) edge (3);
		\path (5) edge (4);
		\path (6) edge (4);
		\end{scope}
		
		\begin{scope}[xshift= 8cm]
		\node (1) at (0, 0) {};
		\node (2) at (0, 1.5) {};
		\node (3) at (0, 3) {};
		\node (4) at (1.5, 0) {};
		\node (5) at (1.5, 1.5) {};
		\node (6) at (1.5, 3) {};
		
		\path (1) edge (2);
		\path (3) edge (2);
		\path (4) edge (5);
		\path (5) edge (6);
		\path (2) edge (5);
		\path (3) edge (6);
		\end{scope}
		
		\begin{scope}[xshift= 12cm]
		\node (1) at (0.75, 0) {};
		\node (2) at (0, 1) {};
		\node (3) at (0.75, 1) {};
		\node (4) at (1.5, 1) {};
		\node (5) at (0.75, 2) {};
		\node (6) at (0.75, 3) {};
		
		\path (1) edge (2);
		\path (3) edge (2);
		\path (4) edge (3);
		\path (5) edge (4);
		\path (6) edge (5);
		\path (2) edge (5);
		\path (1) edge (4);
		\end{scope}
		
		% 		\begin{scope}[yshift= -5cm]
		% 		\node (1) at (0.75, 0) {};
		% 		\node (2) at (0, 1) {};
		% 		\node (3) at (1.5, 1) {};
		% 		\node (4) at (0, 2) {};
		% 		\node (5) at (1.5, 2) {};
		% 		\node (6) at (0.75, 3) {};
		% 		
		% 		\path (1) edge (2);
		% 		\path (3) edge (2);
		% 		\path (4) edge (5);
		% 		\path (5) edge (6);
		% 		\path (4) edge (6);
		% 		\path (1) edge (3);
		% 		\path (2) edge (4);
		% 		\path (5) edge (3);
		% 		\end{scope}
		
		\begin{scope}[xshift= 16cm]
		\node (1) at (0.75, 0) {};
		\node (2) at (0, 1) {};
		\node (3) at (0.75, 1) {};
		\node (4) at (1.5, 1) {};
		\node (5) at (0.75, 2) {};
		\node (6) at (0.75, 3) {};
		
		\path (1) edge (2);
		\path (3) edge (2);
		\path (4) edge (3);
		\path (6) edge (4);
		\path (6) edge (5);
		\path (2) edge (6);
		\path (1) edge (4);
		\path (3) edge (5);
		\end{scope}
		
		% 		\begin{scope}[xshift= 8cm, yshift= -5cm]
		%		\node (1) at (0, 0) {};
		%		\node (2) at (0, 1.5) {};
		%		\node (3) at (0, 3) {};
		%		\node (4) at (1.5, 0) {};
		%		\node (5) at (1.5, 1.5) {};
		%		\node (6) at (1.5, 3) {};
		%		
		%		\path (1) edge (2);
		%		\path (3) edge (2);
		%		\path (4) edge (5);
		%		\path (5) edge (6);
		%		\path (2) edge (5);
		%		\path (3) edge (6);
		%		\path (1) edge (4);
		%		\path (1) edge[bend left] (3);
		%		\path (6) edge[bend left] (4);
		% 		\end{scope}
		
		\begin{scope}[xshift= 20cm]
		\node (1) at (0, 0) {};
		\node (2) at (0, 1.5) {};
		\node (3) at (0, 3) {};
		\node (4) at (1.5, 0) {};
		\node (5) at (1.5, 1.5) {};
		\node (6) at (1.5, 3) {};
		
		\path (1) edge (5);
		\path (3) edge (5);
		\path (4) edge (2);
		\path (6) edge (1);
		\path (2) edge (5);
		\path (3) edge (6);
		\path (1) edge (4);
		\path (2) edge (6);
		\path (3) edge (4);
		\end{scope}
		
		\end{tikzpicture}
		\caption{The six graphs contained in $\cF$. From left to right, we denote them by $F_1, \ldots, F_6$.}
		\label{fig:setF}
	\end{center}
\end{figure}

First, define the set $\cF$ of six bipartite graphs,  see Fig.~\ref{fig:setF}. The smallest one of them is $2P_3$, the largest is $K_{3,3}$, and all the remaining four members of $\cF$ are sandwiched between them. 

%%%%%%%%%%%%%%%%%%%%%%%%%%%

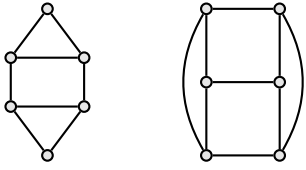
\begin{figure}[!ht]
	\begin{center}
		\begin{tikzpicture}[thick,scale=0.65]
		
		% Define style for nodes
		\tikzstyle{every node}=[circle, draw, fill=black!10,
		inner sep=0pt, minimum width=4pt]
		
		\begin{scope}
		\node (1) at (0.75, 0) {};
		\node (2) at (0, 1) {};
		\node (3) at (1.5, 1) {};
		\node (4) at (0, 2) {};
		\node (5) at (1.5, 2) {};
		\node (6) at (0.75, 3) {};
		
		\path (1) edge (2);
		\path (3) edge (2);
		\path (4) edge (5);
		\path (5) edge (6);
		\path (4) edge (6);
		\path (1) edge (3);
		\path (2) edge (4);
		\path (5) edge (3);
		\end{scope}
		
		\begin{scope}[xshift= 4cm]
		\node (1) at (0, 0) {};
		\node (2) at (0, 1.5) {};
		\node (3) at (0, 3) {};
		\node (4) at (1.5, 0) {};
		\node (5) at (1.5, 1.5) {};
		\node (6) at (1.5, 3) {};
		
		\path (1) edge (2);
		\path (3) edge (2);
		\path (4) edge (5);
		\path (5) edge (6);
		\path (2) edge (5);
		\path (3) edge (6);
		\path (1) edge (4);
		\path (1) edge[bend left] (3);
		\path (6) edge[bend left] (4);
		\end{scope}
		
		\end{tikzpicture}
		\caption{Co-domino (left) and the complement of $C_6$ (right).}
		\label{fig:co-domino}
	\end{center}
\end{figure}

Second, recall that the co-domino graph is the graph shown in Fig.~\ref{fig:co-domino} left, and that the complements of cycles $\overline{C_n}$, $n\ge 5$, are known as  anti-holes; the anti-hole $\overline{C_6}$ is drawn in Fig.~\ref{fig:co-domino} right. 

%%%%%%%%%%%%%%%%%%%%%%%%%%%%%%%%%%%%%%%%%%%%%%%%

\begin{proposition}
	\label{prop:min-imperfect}
	The following graphs are minimally $\gamma_g$-imperfect:
\begin{itemize}
	\item[$(i)$] the path $P_5$;
	\item[$(ii)$] the co-domino; 
	\item[$(iii)$] the anti-hole $\overline{C_n}$ for every $n \ge 5$;
	\item[$(iv)$] each graph from $\cF$. 
	\end{itemize}
\end{proposition}
\proof The graphs referred to in $(i)$, $(ii)$, and $(iv)$ can be checked one-by-one by using Proposition~\ref{lemma:1-2}. For an anti-hole  $\overline{C_n}$ with $n\ge 5$, we first observe that $\gamma(\overline{C_n})=2$ as any two independent vertices form a dominating set. On the other hand, after playing an arbitrary vertex as a first move in the D-game, two non-twin vertices remain undominated and, as follows from  Proposition~\ref{lemma:1-2} $(ii)$, we have  $\gamma_g(\overline{C_n})=3$. Taking an arbitrary proper induced subgraph $H$ of $\overline{C_n}$, either $\gamma(H)=1$ and Proposition~\ref{lemma:1-2} $(i)$ implies $\gamma_g(H)=1$, or there is a vertex with $d_H(v)=|V(H)|-2$ and based on  Proposition~\ref{lemma:1-2} $(iii)$ we may infer $\gamma(H)=\gamma_g(H)=2$. Note that the first case occurs when $\overline{H}$, which is a proper induced subgraph of the cycle $C_n$, contains an isolated vertex; the second case occurs if the minimum degree of $\overline{H}$ equals $1$. Therefore, $\overline{C_n}$ is not a $\gamma_g$-graph but its every proper induced subgraph is a $\gamma_g$-graph. This completes the proof for $(iii)$.
\qed

In the next subsection we will prove that every minimally $\gamma_g$-imperfect graph has $\gamma(G)=2$. We will often refer there to the following statement.
\begin{proposition}
	\label{prop:minimperfect}
If there exist two different vertices $u$ and $v$ in a graph $G$ such that both $N_G\left[u \right] \setminus N_G\left[v \right]$ and $N_G\left[v \right] \setminus N_G\left[u \right]$ contain two nonadjacent vertices, then $G$ is not $2$-$\gamma_g$-perfect.
\end{proposition}
\proof  Suppose  that  $x_1$ and $x_2$ are two independent vertices from $N_G\left[u \right] \setminus N_G\left[v \right]$ and that $y_1$, $y_2$ are two independent vertices from $N_G\left[v \right] \setminus N_G\left[u \right]$. These assumptions directly imply that both $\{v,x_1,x_2\}$ and  $\{u,y_1,y_2\}$ are independent vertex sets in $G$. Hence, the six vertices induce a bipartite graph $H$ with partite classes of size $3$. Checking all the possibilities, we get that $H$ is either isomorphic to a minimally $\gamma_g$-imperfect graph from $\cF$, or contains $P_5$. We may conclude that, under the given conditions, $G$ cannot be $2$-$\gamma_g$-perfect.
\qed

\bigskip

%%%%%%%%%%%%%%%%%%%%%%%%%%%%%%%%%%%%%%%%%%%%%%%%%%%%%
\subsection{Characterization}
\label{subsec:char}
%%%%%%%%%%%%%%%%%%%%%%%%%%%%%%%%%%%%%%%%%%%%%%%%%%%%%

Our goal here is to prove two main results, namely Theorems~\ref{thm:char-1} and \ref{thm:minimp}. Their proofs will be given at the end of the subsection.

\begin{theorem} 
	\label{thm:char-1}
	The following statements are equivalent:
	\begin{itemize}
		\item[$(i)$] $G$ is $\gamma_g$-perfect.
		\item[$(ii)$] $G$ is $2$-$\gamma_g$-perfect.
		\item[$(iii)$] $G$ can be obtained from an isolated vertex by repeatedly applying the following operators:
		\begin{itemize}
			\item For a graph $F$, and for an $s\in \mathbb{N}$, take $F\bigcupdot  K_s$;
			\item For a graph $F$, and for a PSC $\cQ$, take $\cO(F,v,\cQ)$.
		\end{itemize} 
	\end{itemize}
	\end{theorem}

\begin{theorem}
	\label{thm:minimp}
Every minimally $\gamma_g$-imperfect graph has domination number $2$.  
\end{theorem}

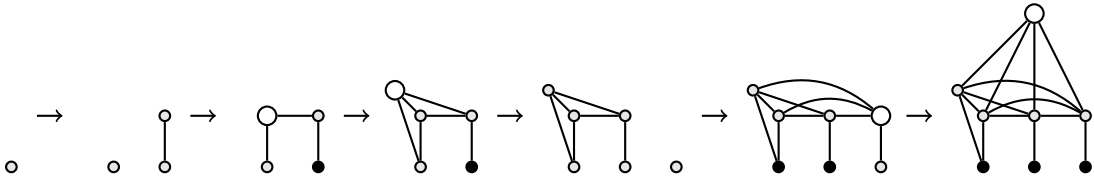
\begin{figure}[!ht]
	\begin{center}
		\begin{tikzpicture}[thick,scale=0.68]
		
		% Define style for nodes
		\tikzstyle{every node}=[circle, draw, fill=black!10,
		inner sep=0pt, minimum width=4pt]
		
		\begin{scope}
		\node (1) at (0,0) {};
		
		\draw [->] (0.5,1) -- (1,1);
		\end{scope}
		
		\begin{scope}[xshift=2cm]
		\node (1) at (0,0) {};
		\node (2) at (1,0) {};
		\node (5) at (1,1) {};
		
		\path (2) edge (5);
		
		\draw [->] (1.5,1) -- (2,1);
		\end{scope}
		
		\begin{scope}[xshift=5cm]
		\node (1) at (0,0) {};
		\node[fill=black] (2) at (1,0) {};
		\node[fill=white, minimum width = 7pt] (4) at (0,1) {};
		\node (5) at (1,1) {};
		
		\path (1) edge (4);
		\path (2) edge (5);
		\path (4) edge (5);
		
		\draw [->] (1.5,1) -- (2,1);
		\end{scope}
		
		\begin{scope}[xshift=8cm]
		\node (1) at (0,0) {};
		\node[fill=black] (2) at (1,0) {};
		\node (4) at (0,1) {};
		\node (5) at (1,1) {};
		\node[fill=white, minimum width = 7pt] (7) at (-0.5,1.5) {};
		
		\path (1) edge (4);
		\path (2) edge (5);
		\path (4) edge (5);
		\path (1) edge (7);
		\path (7) edge (4);
		\path (7) edge (5);
		
		\draw [->] (1.5,1) -- (2,1);
		\end{scope}
		
		\begin{scope}[xshift=11cm]
		\node (1) at (0,0) {};
		\node (2) at (1,0) {};
		\node (3) at (2,0) {};
		\node (4) at (0,1) {};
		\node (5) at (1,1) {};
		\node (7) at (-0.5,1.5) {};
		
		\path (1) edge (4);
		\path (2) edge (5);
		\path (4) edge (5);
		\path (1) edge (7);
		\path (7) edge (4);
		\path (7) edge (5);
		
		\draw [->] (2.5,1) -- (3,1);
		\end{scope}
		
		\begin{scope}[xshift=15cm]
		\node[fill=black] (1) at (0,0) {};
		\node[fill=black] (2) at (1,0) {};
		\node (3) at (2,0) {};
		\node (4) at (0,1) {};
		\node (5) at (1,1) {};
		\node[fill=white, minimum width = 7pt] (6) at (2,1) {};
		\node (7) at (-0.5,1.5) {};
		
		\path (1) edge (4);
		\path (2) edge (5);
		\path (3) edge (6);
		\path (4) edge (5);
		\path (5) edge (6);
		\path (1) edge (7);
		\path (7) edge (4);
		\path (7) edge (5);
		\path (7) edge[bend left] (6);
		\path (4) edge[bend left] (6);
		
		\draw [->] (2.5,1) -- (3,1);
		\end{scope}
		
		\begin{scope}[xshift= 19cm]
		\node[fill=black] (1) at (0,0) {};
		\node[fill=black] (2) at (1,0) {};
		\node[fill=black] (3) at (2,0) {};
		\node (4) at (0,1) {};
		\node (5) at (1,1) {};
		\node (6) at (2,1) {};
		\node (7) at (-0.5,1.5) {};
		\node[fill=white, minimum width = 7pt] (8) at (1,3) {};
		
		\path (1) edge (4);
		\path (2) edge (5);
		\path (3) edge (6);
		\path (4) edge (5);
		\path (5) edge (6);
		\path (1) edge (7);
		\path (7) edge (4);
		\path (7) edge (5);
		\path (7) edge[bend left] (6);
		\path (8) edge (4);
		\path (8) edge (5);
		\path (8) edge (6);
		\path (8) edge (7);
		\path (4) edge[bend left] (6);
		\end{scope}
		
		\end{tikzpicture}
		\caption{A construction of a $\gamma_g$-perfect graph. At the steps where operator $\cO$ is applied, the PSC is marked in black and the newly added vertex is slightly larger.}
		\label{fig:construction}
	\end{center}
\end{figure}

Before we proceed with the proofs, let us consider an example. Figure~\ref{fig:construction} presents a construction of a $\gamma_g$-perfect graph on eight vertices with the consecutive application of operators described in Theorem~\ref{thm:char-1} (iii). Note that by considering the same PSC and just iteratively applying operator $\cO$, we get an infinite family of $\gamma_g$-perfect graphs.

First, we prove that the operator $\cO(G, v, \cQ)$ results in a $\gamma_g$-graph if $G$ is $2$-$\gamma_g$-perfect. We remark that the set of undominated homogeneous cliques referred to at the end of the theorem is not necessarily the same as $\cQ$.

\begin{theorem}
	\label{thm:oO-perfect}
	If $G$ is a $2$-$\gamma_g$-perfect graph and $\cQ$ is a perfect set of cliques in $G$, then  $\cO(G, v, \cQ)$  is a $\gamma_g$-graph. Further, there is an optimal start vertex for Dominator in the D-game on $\cO(G, v, \cQ)$  such that after this first move only a set of homogeneous cliques remains undominated.
\end{theorem}
\proof
Consider graphs $G$ and $G'=\cO(G,v,\cQ)$ that satisfy the conditions of the theorem. First suppose that $\cQ=\emptyset$. Then, by Proposition~\ref{lemma:1-2} (i), we have  $\gamma(G')=\gamma_g(G')$ and the optimal start vertex for Dominator is $v$.
If $\cQ=\{Q_1\}$, then $\gamma(G')$ equals either $1$ or $2$ and, by Proposition~\ref{lemma:1-2} (i) or (ii), the equality $\gamma(G')=\gamma_g(G')$ is true. If $\gamma(G')=1$, there is a universal vertex which is an optimal start vertex; if $\gamma(G')=2$, $v$ is an optimal start vertex and $N_{G'}[v]$ omits only $Q_1$.
\medskip

From now on, we assume that $\cQ=\{Q_1, \ldots, Q_k\}$, where $k\ge 2$. For every $i\in \left[k \right] $ we introduce the notation $N_i'=N'_{G}(Q_i) $ and $N'=\bigcup_{i=1}^k N_i'$. A part of the second neighborhood of a homogeneous clique $Q_i$ is specified as $N_i''=N'_{G}(N_i')\setminus (N'\cup V(\cQ))$; we also define $N''=\bigcup_{i=1}^kN_i''$. By definition, and since any two different cliques $Q_i$ and $Q_j$ are at distance $3$, the sets $N_1',\dots ,N_k'$ are pairwise disjoint and the same is true for the sets $Q_1\cup N_1',\dots, Q_k\cup N_k'$ but it does not hold necessarily for $N_1'', \dots , N_k''$.  Observe further that, by definition of PSC and by the condition $k\ge 2$, the set $N_i'$ is not empty for any  $i \in \left[ k \right] $.

 By Proposition~\ref{prop:opO}, $\gamma(G')=\gamma_g(G')$ is true if $\gamma(G) > k$. Note that in this case $v$ is an optimal start vertex in the D-game and $V(G')\setminus N_{G'}[v]=V(\cQ)$.
 \medskip
 
 From now on, we may assume that there is a dominating set $D$ of cardinality $k$ in $G$. We will also suppose that, under this condition, $D$ is chosen such that $|D\cap N'|$ is maximum. To dominate all the cliques $Q_1, \dots, Q_k$, the set $D$ must contain at least one vertex from each  $Q_i\cup N_i'$. Consequently, $D \subseteq V(\cQ) \cup N' $ and $D$ contains exactly one vertex, say $d_i$, from each $Q_i\cup N_i'$. 

Under the condition $\gamma(G)=k \ge 2$, we continue the proof with a series of claims.

\begin{claim}
	\label{claim:1}
	There exists a vertex $y$ in $N'\cap D$ such that $N'' \subseteq N[y]$.
\end{claim}
\proof  
Suppose that $d_i$ and $d_j$ are different vertices with $d_i\in N_i'\cap D$, and $d_j \in N_j' \cap D$ and let $x_i$ be a vertex from $Q_i$. If $d_i$ has a neighbor $z_i\in N''$ which is not a neighbor of $d_j$, then $x_i$ and $z_i$ are two independent vertices from $ N\left[ d_i\right] \setminus  N\left[ d_j\right]$. Similarly, if $d_j$ has a neighbor $z_j$ in $N''$ which is not a neighbor of $d_i$, then $ N\left[ d_j\right] \setminus  N\left[ d_i\right]$ contains two independent vertices. These two facts together, by Proposition~\ref{prop:minimperfect}, would contradict the $2$-$\gamma_g$-perfectness of $G$. Hence, at least one of $(N\left[ d_i\right] \setminus N\left[ d_j\right])\cap N''$ and $(N\left[ d_j\right] \setminus N\left[ d_i\right])\cap N''$ is empty, hence $N\left[ d_i\right]\cap N''$ is a subset of $N\left[ d_j\right]\cap N''$ or vice versa. This defines a linear ordering for the sets $N\left[ d_i\right]\cap N''$, $i\in [ k]$, even if $|N'\cap D| = 1$. As $D$ dominates all vertices from $N''$,  there exists a vertex in $N'\cap D$ that dominates the entire $N''$. \smallqed

\medskip

\begin{claim}
	\label{claim:2}
	There is at most one $i\in [k]$ such that the set $N'_i$ does not induce a complete subgraph in $G$. 
\end{claim}
\proof Suppose for a contradiction that neither $G[N_i']$ nor $G[N_j']$ is complete. Then, for every $x_i\in Q_i$ and $x_j\in Q_j$, both sets $N[x_i]\setminus N[x_j]$ and $N[x_j]\setminus N[x_i]$ contain nonadjacent vertices. By Proposition~\ref{prop:minimperfect}, this contradicts the $2$-$\gamma_g$-perfectness of $G$. \smallqed
\medskip

\begin{claim}
	\label{claim:3}
	If $N''=\emptyset$, then $\gamma(G')=\gamma_g(G')$. 
\end{claim}
\proof By Claim~\ref{claim:2} and since $k \ge 2$, there is a set $N'_j$ that induces a complete subgraph in $G$. As $\cQ$ is a perfect set of cliques in $G$, any vertex $y_j \in N'_j$ dominates the entire $N' \cup Q_j \cup \{v\}$ in $G'$. If Dominator first plays such a vertex $y_j$, only the $k-1$ homogeneous cliques different from $Q_j$ remain undominated because $\gamma(G) = k$. In the continuation of the game, under any strategy of Staller and Dominator, exactly one homogeneous clique will be dominated with each move. Therefore, this (optimal) first move of Dominator ensures that the game finishes within $k$ moves. Consequently, we have $\gamma_g(G') \le k=\gamma(G)$ and  may conclude $\gamma(G')=\gamma_g(G')$. \smallqed

\medskip

\begin{claim}
	\label{claim:4}
	If $N''\neq\emptyset$, then $\gamma(G')=\gamma_g(G')$. 
\end{claim}
\proof 
First, suppose that there is no incomplete $N'_i$ and, consequently, $N'$ induces a complete subgraph in $G'$. Then, by Claim~\ref{claim:1}, Dominator may choose a  vertex $y \in N'$ which dominates the entire $N'' \cup N'\cup\{v\}$ and also dominates one homogeneous clique $Q_i$. Note that, by the condition $\gamma(G)=|\cQ|$, all the vertices of $G'$ are contained in $V(\cQ)\cup N' \cup N'' \cup \{v\}$. Thus, after the first move $y$ of Dominator, only $k-1$ cliques from $\cQ$ remain undominated and the game will be finished with $k$ moves. This proves $\gamma(G')=\gamma_g(G')$.

In the other case, when $N'$ is not complete, we might also have a vertex from $ N'$ that dominates the entire $N'' \cup N'\cup\{v\}$. This implies $\gamma(G')=\gamma_g(G')$, again.

What remains is to consider the case when we do not have a vertex in $N'$ which dominates the entire $N' \cup N''$. By Claim~\ref{claim:1}, we have a vertex $y \in N'$ that dominates $N''$. We may assume, without loss of generality, that $y\in N_1'$. Then there is a vertex $y'\in N_1'$ which is not dominated by $y$. Since by Claim~\ref{claim:2}, $N_2'$ induces a clique, any vertex $z$ from $N_2'$ dominates the entire $N'$ but does not dominate $N''$. Fixing any such vertex $z$, there exists a vertex $w$ from $N''$ which is nonadjacent to $z$ but adjacent to $y$. Further, let $x_1 \in Q_1$ and  $x_2 \in Q_2$ be two arbitrary vertices from the homogeneous cliques. Observe that the two independent vertices $x_1$ and $w$  belong to $N_{G}[y]\setminus N_{G}[z]$  and also that the independent vertices $x_2$ and $y'$ are contained in $N_{G}[z]\setminus N_{G}[y]$. By Proposition~\ref{prop:minimperfect}, this case is not possible as it contradicts the $2$-$\gamma_g$-perfectness of $G$. \smallqed

\medskip

Our previous discussions on the cases $|\cQ|=0$ and $|\cQ|=1$, Proposition~\ref{prop:opO}, Claims~\ref{claim:3} and \ref{claim:4} together imply $\gamma(G')=\gamma_g(G')$ for any $G'= \cO(G,v,\cQ)$, where $G$ is $2$-$\gamma_g$-perfect. An optimal start vertex with the required property was identified for all cases. This finishes the proof of Theorem~\ref{thm:oO-perfect}. 
\qed

After proving a lemma, we will show that every $2$-$\gamma_g$-perfect graph $G$ can be obtained from another 2-$\gamma_g$-perfect graph $F$ by using the operator disjoint union with a complete graph or the operator $\cO$. 
%In fact, here we put a (seemingly) weeker condition, where only the induced subgraphs of domination number $2$ are required to be $\gamma_g$-graphs.
We say that $F$ is a \emph{$\gamma$-$2$-maximal subgraph} of $G$ if it is an induced subgraph of $G$ with $\gamma(F)\le 2$  and inclusion-wise maximal with this property. That is, for any induced subgraph $F'$ of $G$ with $V(F) \subsetneqq V(F')$ we have $\gamma(F')\ge 3$. 

\begin{lemma}
	\label{lemma:2-max}	
	Let $G$ be a  $2$-$\gamma_g$-perfect graph and let $v$ be a vertex in $G$ such that $d_{\widehat{G}}(\widehat{v})=\Delta(\widehat{G})$. If $F$ is a $\gamma$-$2$-maximal subgraph of $G$ which contains the entire $N[v]$, then $\widehat{v}$ is a vertex of maximum degree in $\widehat{F}$ and $v$ is an optimal start vertex in the D-game on $F$.
\end{lemma}

\proof Under the given conditions, $\gamma(F)=\gamma_g(F)=2$. By Observation~\ref{observation} and Proposition~\ref{lemma:1-2}, this implies $\Delta(\widehat{F})=|V(\widehat{F})|-2$ and hence, there is a vertex $\widehat{u}$ which is adjacent to all but one vertex of $\widehat{F}$. Now, assume for a contradiction that $d_{\widehat{F}}(\widehat{v}) < d_{\widehat{F}}(\widehat{u}) $. We consider two cases under this assumption.
\begin{itemize}
	\item First suppose that $N_G[v] \subseteq N_G[u]$. This implies $N_{\widehat{G}}[\widehat{v}] \subseteq N_{\widehat{G}}[\widehat{u}]$ and, since $\widehat{v}$ is of maximum degree in $\widehat{G}$, we have $d_{\widehat{G}}(\widehat{v}) = d_{\widehat{G}}(\widehat{u})$. Consequently, $N_{\widehat{G}}[\widehat{v}] = N_{\widehat{G}}[\widehat{u}]$ and $N_G[v]=N_G[u]$ hold. The latter equality implies $N_F[v]=N_F[u]$ that contradicts the assumption $d_{\widehat{F}}(\widehat{v}) < d_{\widehat{F}}(\widehat{u}) $.
	\item In the other case, we suppose that $u$ is not adjacent to all vertices of $N_G[v]$. Since $\widehat{u}$ has degree $|V(\widehat{F})|-2$ in $\widehat{F}$, there is exactly one maximal homogeneous clique $Q$ in $F$ such that $Q \subseteq N_G[v]\setminus N_G[u]$. Note that, since $\gamma(F)=2$, there exist some vertices in $F$ (and in $G$) which are adjacent to $u$ but not adjacent to $v$. Then, $d_{\widehat{F}}(\widehat{v}) < d_{\widehat{F}}(\widehat{u}) $ and $d_{\widehat{G}}(\widehat{v}) \ge d_{\widehat{G}}(\widehat{u}) $ imply that the  clique $Q$ is homogeneous in $F$ and non-homogeneous in $G$. That is, there are two vertices, say $x_1$ and $x_2$ in $Q$  such that $x_1$ has a neighbor $z$ with $zx_2 \notin E(G)$ and $z \notin V(F)$. Then, $N_G[u]\cup N_G[x_1]$ contains $V(F)$ as a proper subset. The subgraph $G[N_G[u]\cup N_G[x_1]]$ is dominated by $u$ and $x_1$ and, therefore, $F$ cannot be a  $\gamma$-$2$-maximal subgraph of $G$. 
\end{itemize} 
 As both possible cases were concluded with contradictions, we infer that $\widehat{v}$ is a vertex of maximum degree in $\widehat{F}$ and, since $\widehat{F}$ is a $\gamma_g$-graph, $d_{\widehat{F}}(\widehat{v})= |V(\widehat{F})|-2$. Thus, in the D-game on $F$, the vertex $v$ is an optimal start vertex for Dominator. \qed

\begin{theorem}
	\label{thm:construction-perf}
	Every  $2$-$\gamma_g$-perfect graph $G$ can be constructed in the following way:
	\begin{itemize}
		\item[$(i)$] If $G$ is disconnected, then it can be obtained as $G' \bigcupdot  K_s$, where $G'$ is also $2$-$\gamma_g$-perfect.
		\item[$(ii)$] If $G$ is connected, then for every vertex $v\in V(G)$ with $d_{\widehat{G}}(\widehat{v})=\Delta(\widehat{G})$, there exists a perfect set of cliques $\cQ$ in $G-v$ so that $G=\cO(G-v, v, \cQ)$. Further, $G-v$ is a $2$-$\gamma_g$-perfect graph.
			\end{itemize} 
\end{theorem}

\proof $(i)$ Since every induced subgraph $F$ of $G$ of domination number $2$ is a $\gamma_g$-graph, $G$ is $2P_3$-free and hence, it cannot contain more than one non-complete component. This yields that a disconnected  $G$ can always be obtained as a disjoint union $G' \bigcupdot  K_s$. Since $G'$ is an induced subgraph of $G$, it is $2$-$\gamma_g$-perfect as well.

$(ii)$ Consider a connected graph $G$ satisfying the conditions in the theorem and consider a vertex $v$ with $d_{\widehat{G}}(\widehat{v})=\Delta(\widehat{G})$. We prove that the vertices outside $N_G[v]$ form a PSC in $G$. It is clearly true, if $\gamma(G)=1$ and $v$ is a universal vertex.

First, suppose that there are two vertices $x$ and $x'$ in $V(G)\setminus N_G[v]$ such that $x$ and $x'$ are not twins in $G$ but they are adjacent. Then, we have a vertex $z$ in $G$ which is adjacent to exactly one of $x$ and $x'$; we may suppose that $xz \in E(G)$ and $x'z \notin E(G)$. Since $\{v,x\}$ is a 2-element dominating set in the subgraph induced by $N_G[v] \cup \{z,x,x'\}$, we may consider a $\gamma$-$2$-maximal subgraph $F$ of $G$ which contains all the vertices from  $N_G[v] \cup \{z,x,x'\}$. By Lemma~\ref{lemma:2-max}, the vertex $\widehat{v}$ must be of maximum degree in $\widehat{F}$ and in particular, $d_{\widehat{F}}(\widehat{v}) = |V(\widehat{F})|-2$ must hold. This contradicts the fact that $v$ is not adjacent in $F$ to at least two non-twin vertices, namely to $x$ and $x'$. This contradiction proves that $V(G)\setminus N_G[v]$ consists of components $Q_1, \dots, Q_k$ which are homogeneous cliques in $G$.

Secondly, we prove that no two of the homogeneous cliques $Q_1, \dots, Q_k$ are at distance $2$. Suppose, to the contrary, that two vertices, $x_i$ from $Q_i$ and $x_j$ from $Q_j$, where $i\neq j$, have a common neighbor $y$ in $G$. All vertices in $N_G[v] \cup \{x_i, x_j\}$ are dominated by $v$ and $y$. Therefore, we may consider again a $\gamma$-$2$-maximal subgraph $F$ of $G$ which contains all the vertices from  $N_G[v] \cup \{x_i,x_j\}$. By Lemma~\ref{lemma:2-max},  $d_{\widehat{F}}(\widehat{v}) = |V(\widehat{F})|-2$ must hold. On the other hand, as the non-twin vertices $x_i$ and $x_j$ are outside $N_G[v]$, we have $d_{\widehat{F}}(\widehat{v}) \le |V(\widehat{F})|-3$. This contradiction proves that $d(Q_i, Q_j) \ge 3$ for any two different indices $i$ and $j$. This is true when either $G$ or $G-v$ is considered.

Finally, observe that there is an edge between any two vertices $y_i$ and $y_j$ whenever  $y_i \in N_i'$ and $y_j \in N_j'$, $i\neq j$. Indeed, in case of  $y_iy_j \notin E(G)$ we would have an induced $P_5$, namely  $x_iy_ivy_jx_j$ in $G$. Since $\gamma(P_5)=2 <\gamma_g(P_5)$, this contradicts the condition in the theorem. We conclude that, under the conditions of part $(ii)$,  $\cQ=\{Q_1,\dots , Q_k\}$ is a PSC in $G-v$ and $G$ can be obtained as $\cO(G-v, v, \cQ)$. \qed

\begin{proposition}
	\label{prop:disj-union}
If $G$ is a $2$-$\gamma_g$-perfect graph, then $G\bigcupdot  K_s$ is a  $\gamma_g$-graph.
\end{proposition}
\proof If $G$ consists of $c$ complete components, then $\gamma(G\bigcupdot  K_s)=\gamma_g(G\bigcupdot  K_s)=c+1$. Otherwise, $G$ contains $c$ complete components ($c\ge 0$) and  exactly one component, say  $G'$, which is not a complete graph. Since $G'$ is $2$-$\gamma_g$-perfect and connected, by Theorem~\ref{thm:construction-perf} $(ii)$ and Theorem~\ref{thm:oO-perfect}, there is an optimal first move $v$ for Dominator such that $V(G')\setminus N[v]$ consists of $k$ homogeneous cliques. If Dominator plays this vertex  $v$ as his first move on the entire $G\bigcupdot  K_s$, then $k+c+1$ homogeneous cliques remain undominated and we have $\gamma(G\bigcupdot  K_s)=\gamma_g(G\bigcupdot  K_s)= k+c+2$. This proves that $G\bigcupdot  K_s$ is a $\gamma_g$-graph. \qed

%\begin{theorem}
%	\label{thm:perf-2-perf}
%	$G$ is a $2$-$\gamma_g$-perfect graph if and only if it is  $\gamma_g$-perfect.
%\end{theorem}
%\proof

%Let $v \in V(G)$ be a vertex with $d_{\widehat{G}}(\widehat{v})=\Delta(\widehat{G})$. If $G$ satisfies the condition on the induced subgraphs with domination number $2$, then the same condition holds for the subgraph $F=G-v$. By the induction hypothesis, $F$ is $\gamma_g$-perfect and, by part $(i)$ of the present theorem, $G$ can be obtained as $\cO(F, v, \cQ)$. Then, Theorem~\ref{thm:oO-perfect} implies that $G$ is $\gamma_g$-perfect. \qed

Now we are ready to prove the characterization theorem.
\medskip

\noindent \textbf{Proof of Theorem~\ref{thm:char-1}.} First, we show that $(ii)$ implies $(i)$. Note that the statement is clearly true for graphs of small order and then, we may proceed by induction on the number of vertices in $G$.  Suppose that $G$ is $2$-$\gamma_g$-perfect. By definition, it is also true for every induced subgraph of $G$. Hence, for every proper induced subgraph $F$ of $G$, the induction hypothesis implies that $F$ is $\gamma_g$-perfect. As for $G$ itself, if it is disconnected then, by Theorem~\ref{thm:construction-perf} $(i)$, $G$ can be obtained as $G' \bigcupdot  K_s$ from a $2$-$\gamma_g$-perfect $G'$.  Proposition~\ref{prop:disj-union} then implies $\gamma(G)=\gamma_g(G)$. In the other case, $G$ is connected and, by Theorem~\ref{thm:construction-perf} $(ii)$ and Theorem~\ref{thm:oO-perfect}, we have $\gamma(G)=\gamma_g(G)$ again. This proves that every $2$-$\gamma_g$-perfect graph is $\gamma_g$-perfect, that is, $(ii) \Rightarrow (i)$. The other direction immediately follows from the definitions. Therefore, $G$ is $\gamma_g$-perfect if and only if it is $2$-$\gamma_g$-perfect.

Now, we prove that $(ii)$ is equivalent with $(iii)$. By Theorem~\ref{thm:construction-perf}, each $2$-$\gamma_g$-perfect graph $G$ can be obtained from an appropriate smaller $2$-$\gamma_g$-perfect graph $G'$ as $G=G' \bigcupdot K_s$ or as $G=\cO(G', v , \cQ)$. Applying the theorem for the $2$-$\gamma_g$-perfect $G'$ and then, repeatedly, for the smaller graphs, the process ends with $K_1$. This proves the  implication $(ii) \Rightarrow (iii)$. 

To prove the other direction, we proceed by induction on the order of the graph. Suppose that a graph $H$ can be built from $K_1$ by using the two operators specified in $(iii)$. As $K_1$ is a $\gamma_g$-graph and Theorem~\ref{thm:oO-perfect} and Proposition~\ref{prop:disj-union} say that the operators preserve this property, $H$ is a $\gamma_g$-graph. We consider an arbitrary induced proper subgraph $F$ of $H$ which satisfies $\gamma(F)=2$.
\begin{itemize}
\item If $H=H' \bigcupdot K_s$ and $H'$ can be built from $K_1$ by using the two operators, then, by the induction hypothesis, $H'$ is $2$-$\gamma_g$-perfect. If $F$ is an induced subgraph of $H'$,  then $F$ is also $2$-$\gamma_g$-perfect. In the other case, since $\gamma(F)=2$, $F$ contains some vertices from $K_s$ and meets $H'$ in a subgraph of domination number $1$. In either case, $F$ is a $\gamma_g$-graph.
\item Suppose that $H=\cO(H',v, \cQ)$, where $H'$ is built from $K_1$ by using the two operators.  We have two cases again. First, assume that $F$ is a subgraph of $H'$. Since $H'$ is $2$-$\gamma_g$-perfect by the induction hypothesis, the same is true for $F$. In the second case, $F$ contains $v$ and also contains some (outer) vertices which are outside $N_H[v]$, because otherwise it would hold $\gamma(F) = 1$. Since $\gamma(F)=2$, these outer vertices belong to either one or two homogeneous cliques from $\cQ$. If $F$ meets only one homogeneous clique from $\cQ$, then, no matter whether $F$ is connected or not, $v$ is an optimal start vertex in the D-game on $F$ and we have $\gamma_g(F)=\gamma(F)=2$. Finally, suppose that $F$ meets two cliques from $\cQ$, say $Q_i$ and $Q_j$. If $F$ contains vertices from both $N'_i=N'_H(Q_i)$ and $N'_j=N'_H(Q_j)$, then $F$ also contains the join between $N_i'\cap F$ and $N_j'\cap F$. Thus, $V(F)\cap V(H')$ induces a subgraph $H''$ of $H'$ such that $\cQ'=\{V(F) \cap Q_i, V(F) \cap Q_j\}$ is a PSC in $H''$ and  $F=\cO(H'', v, \cQ')$. Since $H''$ is an induced subgraph of $H'$, it is $2$-$\gamma_g$-perfect, and  Theorem~\ref{thm:oO-perfect}  implies that $F$ is a $\gamma_g$-graph. If $F$ does not meet $N_i'$ but, as it was assumed, $F$ meets $Q_i$, then $F$ is disconnected. In this case, $V(F) \cap Q_i$ induces a clique component in $F$, while the another component must have a universal vertex. This universal vertex is an optimal start vertex in the D-game and we have $\gamma_g(F)=\gamma(F)=2$.
\end{itemize}
Hence, in $H$, every induced subgraph $F$ of domination number $2$ is a $\gamma_g$-graph. In other words, $H$ is a $2$-$\gamma_g$-perfect graph. We may conclude that $(iii)$ implies $(ii)$ and, therefore, $(i)$, $(ii)$, and $(iii)$ are equivalent. \qed

\noindent \textbf{Proof of Theorem~\ref{thm:minimp}.} Assume for a contradiction that there exists a graph $G$ with $\gamma(G)\ge 3$ which is minimally $\gamma_g$-imperfect. That is, every proper induced subgraph of $G$, including all induced subgraphs of domination number $2$, are $\gamma_g$-perfect. But then, Theorem~\ref{thm:char-1} implies that $G$ is also $\gamma_g$-perfect. This contradiction completes the proof. \qed

%%%%%%%%%%%%%%%%%%%%%%%%%%%%%%%%%%%%%%%%%%%%%%%%%%%
%%%%%%%%%%%%%%%%%%%%%%%%%%%%%%%%%%%%%%%%%%%%%%%%%%%%%
\section{Some applications}
\label{sec:appl}
%%%%%%%%%%%%%%%%%%%%%%%%%%%%%%%%%%%%%%%%%%%%%%%%%%%%%
%%%%%%%%%%%%%%%%%%%%%%%%%%%%%%%%%%%%%%%%%%%%%%%%%%%

In this section we present some applications of the characterization from Section~\ref{sec:gamma-g}.

%%%%%%%%%%%%%%%%%%%%%%%%%%%%%%%%%%%%%%%%%%%%%%%%%%%%%
\subsection{Recognition complexity}
\label{subsec:recognition}
%%%%%%%%%%%%%%%%%%%%%%%%%%%%%%%%%%%%%%%%%%%%%%%%%%%%%

Here we demonstrate that the characterizations of $\gamma_g$-perfect graphs are constructive, more precisely we have: 

\begin{theorem}
\label{thm:recognition}
Graphs that are $\gamma_g$-perfect can be recognized in polynomial time.
\end{theorem}

\proof
Let $G$ be an arbitrary graph. Then we first determine its connected components and in view of Theorem~\ref{thm:char-1} discard the components that induce complete graphs. Clearly, this can be done in polynomial time. Then $G$ is $\gamma_g$-perfect if and only if the remaining non-complete component $H$ is $\gamma_g$-perfect. Let $v\in V(H)$ be an arbitrary vertex with the property that $d_{\widehat{H}}(\widehat{v})=\Delta(\widehat{H})$. Clearly, such a vertex can be found in polynomial time. By Theorem~\ref{thm:construction-perf}, $H$ is $2$-$\gamma_g$-perfect (and hence $\gamma_g$-perfect by Theorem~\ref{thm:char-1}) if and only if $H-N[v]$ consists of a perfect set of cliques in $H-v$ and $H-v$ is $2$-$\gamma_g$-perfect. Since checking whether $H-N[v]$ consists of a perfect set of cliques can be done in polynomial time, we have reduced in polynomial time the problem of whether $H$ is $\gamma_g$-perfect to verifying whether $H-v$ is $\gamma_g$-perfect. Repeating the procedure on $H-v$ yields a polynomial recognition algorithm. 
\qed

%%%%%%%%%%%%%%%%%%%%%%%%%%%%%%%%%%%%%%%%%%%%%%%%%%%%%
\subsection{Triangle-free graphs}
\label{subsec:triangle-free}
%%%%%%%%%%%%%%%%%%%%%%%%%%%%%%%%%%%%%%%%%%%%%%%%%%%%%

In this section, we study triangle-free $\gamma_g$-perfect and minimally $\gamma_g$-imperfect graphs. In particular, we also discuss trees. For this, we need the following notation. A graph $KC_{m,n}$, $m \geq 1, n \geq 0$, has vertices $\{c, d, u_1, \ldots, u_m, v_1, \ldots, v_n \}$ and edges $c \sim u_j \sim d$, $j \in [m]$, and $c \sim v_i$, $i \in [n]$. Additionally, set $KC_{0, n} = K_{1,n}$. Note that the graph $KC_{1, n}$ is the star $K_{1,n}$ with a pendant $P_2$ attached to the central vertex, and graphs $KC_{m, n}$, $m \geq 2$, contain 4-cycles.

\begin{proposition}
	\label{prop:triangle-free-perfect}
	A connected triangle-free graph is $\gamma_g$-perfect if and only if it is isomorphic to $KC_{m,n}$ for some $m, n \geq 0$.
\end{proposition}

\proof
It can easily be checked, that graphs $KC_{m,n}$, for some $m, n \geq 0$, are all connected, triangle-free, and $\gamma_g$-perfect. 

Now suppose $G$ is a connected triangle-free $\gamma_g$-perfect graph. By Theorem~\ref{thm:char-1}, $G = \cO(F, v, \cQ)$ for some triangle-free $\gamma_g$-perfect graph $F$. As $G$ is triangle-free, $|\cQ| \in \{0, 1\}$ and each clique in $\cQ$ can have only one vertex. In addition, the vertices $V(F) \setminus V(\cQ)$ must induce an independent set, otherwise $G$ would again contain a triangle. 

If $|\cQ| = 0$, then $F$ is just a graph with no edges, thus $\cO(F, v, \cQ)$ is a star, i.e.\ a graph $KC_{0, n}$. If $|\cQ| = 1$, then $F$ can only be a star with some additional independent vertices, where the only clique in $\cQ$ consists of the center of the star. In this case $\cO(F, v, \cQ)$ is isomorphic to $KC_{m, n}$, $m \geq 1$.
\qed

\begin{corollary}
	\label{cor:trees-perfect}
	A tree $T$ is $\gamma_g$-perfect if and only if it is isomorphic to $KC_{0, n}$ or $KC_{1, n}$.
\end{corollary}

Note that the last corollary can also be obtained from the characterization of $\gamma_g$-minimal trees~\cite{nadjafi-2016}.

In the rest of the subsection we investigate minimally $\gamma_g$-imperfect triangle-free graphs. We first prove a more general result.

\begin{theorem}
	\label{thm:connected-imperfect}
	The only disconnected minimally $\gamma_g$-imperfect graph is $2 P_3$.
\end{theorem}

\proof
Let $G$ be a disconnected minimally $\gamma_g$-imperfect graph. From Theorem~\ref{thm:minimp} it follows that $\gamma(G) = 2$. Thus $G$ has exactly two connected components. At least one of these two components is non-complete, otherwise $\gamma_g(G) = 2$ and $G$ is not imperfect. If both components are non-complete, then either $G = 2 P_3$, or $G$ contains $2 P_3$ as an induced subgraph and is thus not a minimally $\gamma_g$-imperfect. The only remaining case is that one of the components is complete, and the other component, denote it with  $G'$, is not. But then if Dominator starts the domination game by playing the universal vertex of $G'$, Staller has no other option but to finish the game. Thus $\gamma_g(G) = 2$ and $G$ is not imperfect.
\qed

We now list all minimally $\gamma_g$-imperfect trees and then use this result to determine all triangle-free minimally $\gamma_g$-imperfect graphs.

\begin{proposition}
	\label{prop:trees-imperfect}
	The only minimally $\gamma_g$-imperfect trees are $P_5$ and the tree $F_2$ from the family $\cF$.
\end{proposition}

\proof
Let $v$ be a leaf of a minimally $\gamma_g$-imperfect tree $T$. Then $T-v$ is a perfect tree, hence it is by Corollary~\ref{cor:trees-perfect} either $KC_{0, n}$ or $KC_{1, n}$. If $T-v = KC_{0, n}$, then attaching $v$ either to $c$ or to $v_i$, we get a perfect tree. If $T-v = KC_{1, n}$, then attaching $v$ to $c$ yields a perfect tree. On the other hand, attaching $v$ to $d$ or $v_i$ results in a graph, which contains $P_5$ as an induced subgraph. Thus $T$ is minimally $\gamma_g$-imperfect only if it is equal to $P_5$. The remaining case is that $v$ is attached to $u_1$. In this case, $T$ contains the tree $F_2$ from the family $\cF$ as an induced subgraph, hence it is minimally $\gamma_g$-imperfect only if it is isomorphic to this tree.
\qed

\begin{theorem}
	\label{thm:triangle-free-imperfect}
	The only minimally $\gamma_g$-imperfect triangle-free graphs are $P_5$, $C_5$ and the graphs from $\cF$.
\end{theorem}

\proof
Let $G$ be a minimally $\gamma_g$-imperfect triangle-free graph. If $G$ is disconnected, then the only possibility is $F_1 = 2 P_3$ by Theorem~\ref{thm:connected-imperfect}. If all non-pendant vertices of a connected graph $G$ are cut vertices, then $G$ is a tree and Proposition~\ref{prop:trees-imperfect} settles this case. The only remaining case is that $G$ is connected and has a vertex $v$ which is not a cut vertex. This means that $G-v$ is a connected triangle-free $\gamma_g$-perfect graph, thus by Proposition~\ref{prop:triangle-free-perfect}, $G-v = KC_{m,n}$ for some $m,n\geq 0$. We consider different possibilities for $G-v$ and analyze how $v$ can be added back to obtain the graph $G$. Note that as $G$ is triangle-free, the neighborhood of $v$ is an independent set. We use the notation introduced at the beginning of Subsection~\ref{subsec:triangle-free} and additionally set $\cV = \{v_1, \ldots, v_n\}$ and $\cU = \{ u_1, \ldots, u_m \}$. Recall also the graphs from the family $\cF = \{ F_1, \ldots, F_6\}$, see Fig.~\ref{fig:setF} again.

Note that whenever we determine a minimally $\gamma_g$-imperfect subgraph $H$ of $G$, the only possibility for $G$ to be minimally $\gamma_g$-imperfect is that $H = G$. Hence in the remaining part of the proof, we only determine minimally $\gamma_g$-imperfect subgraphs. 

If $G-v = KC_{0,n}$, then, if $v$ is connected only to $c$ or if $v$ is connected to some vertices in $\cV$, we see that $G$ is perfect. If $G-v = KC_{1,n}$, then we only need to consider cases where $v$ belongs to a cycle in $G$ (otherwise $G$ is a tree and Propositon~\ref{prop:trees-imperfect} settles this case). If $v$ is adjacent to $c$ and $d$, then $G$ is perfect. If $N(v) \subseteq \cV$, then $G$ contains an induced $P_5$. If $v$ is adjacent to $d$ and some vertices in $\cV$, then $G$ contains an induced $C_5$. The remaining options yield a graph with triangles.

We now study the case $G-v = KC_{m,n}, m \geq 2$. We distinguish the following cases.

\medskip\noindent
\textbf{Case 1}: $v$ has only one neighbor in $G$.
	\begin{itemize}
		\item $n = 0$:\\
		If $v$ is adjacent to $c$ or $d$, then $G$ is perfect. The same holds if $v$ is adjacent to $u_1$ or $u_2$, and $m = 2$. But if $m \geq 3$ and $v$ is adjacent to one of the vertices in $\cU$, then $G$ contains an induced $F_4$. 
		
		\item $n \geq 1$:\\
		If $v$ is adjacent to $c$, then $G$ is perfect. If $v$ is adjacent to $d$ or a vertex in $\cV$, then $G$ contains an induced $P_5$. If $v$ is adjacent to a vertex in $\cU$, then $G$ contains an induced $F_3$. 
	\end{itemize}
\medskip\noindent
\textbf{Case 2}: $v$ has at least two neighbors in $G$.
	\begin{itemize}
		\item $n = 0$:\\
		If $m = 2$, then all cases can be easily checked, and we see that non of the obtained graphs is minimally $\gamma_g$-imperfect. Hence we suppose $m \geq 3$. 
		
		If $v$ is adjacent to $c$ and $d$, then $G$ is perfect. If $v$ is adjacent to exactly two vertices in $\cU$, then $G$ contains an induced $F_5$. If $v$ is adjacent to three or more vertices in $\cU$, then $G$ contains an induced $F_6$. 
		
		\item $n \geq 1$:\\
		If $v$ is adjacent to $c$ and $d$, then $G$ is perfect. If $N(v) \subseteq \cV \cup \cU$ and $v$ is not adjacent to all vertices in $\cU$, then $G$ contains one of $P_5$, $F_3$, or $F_5$ as an induced subgraph. If $\cU \subseteq N(v) \subseteq \cV \cup \cU$, then there are a few cases. If $m = 2$, then either $G$ contains an induced $F_4$ (if $v$ has no neighbors in $\cV$) or $G$ contains an induced $F_5$ (if $v$ has at least one neighbor in $\cV$). If $m \geq 3$, then $G$ contains an induced $F_6$. And finally, if $v$ is adjacent to $d$ and some vertices in $\cV$, then $G$ contains an induced $C_5$. 
	\end{itemize}

It follows from the above case analysis, that $G$ is minimally $\gamma_g$-imperfect only if it is equal to one of the graphs from $\{P_5, C_5\} \cup \cF$.
\qed

%%%%%%%%%%%%%%%%%%%%%%%%%%%%%%%%%%%%%%%%%%%%%%%%%%%
%%%%%%%%%%%%%%%%%%%%%%%%%%%%%%%%%%%%%%%%%%%%%%%%%%%
\section{Further types of perfectness}
\label{sec:further}
%%%%%%%%%%%%%%%%%%%%%%%%%%%%%%%%%%%%%%%%%%%%%%%%%%%
%%%%%%%%%%%%%%%%%%%%%%%%%%%%%%%%%%%%%%%%%%%%%%%%%%%

\begin{proposition}
	A graph $G$ is $\gamma_g'$-perfect, if and only if it is the disjoint union of cliques.
\end{proposition}
\proof Since $\gamma(P_3)=1 <\gamma_g'(P_3)=2$, $P_3$ is a minimally $\gamma_g'$-imperfect graph. Then, $G$ is $P_3$-free; it is a disjoint union of some, say $k$, clique components. It is easy to see that $\gamma(G)=\gamma_g'(G)=k$. \qed

\begin{proposition}
\label{prp:gamma-tg-perfect}
	An isolate-free graph $G$ is $\gamma_{tg}$-perfect, if and only if it is $(P_4, \overline{2P_3})$-free.
\end{proposition}

\proof 
Since $\gamma_t(P_4) = 2 < 3 = \gamma_{tg}(P_4)$ and $\gamma_t(\overline{2P_3}) = 2 < 3 = \gamma_{tg}(\overline{2P_3})$, it follows that an isolate-free, $\gamma_{tg}$-perfect graphs is $(P_4, \overline{2P_3})$-free.

Suppose now that $G$ is an isolate- and $(P_4, \overline{2P_3})$-free graph. To prove that $G$ is $\gamma_{tg}$-perfect we proceed by induction on $n(G)$, where the cases with $n(G)\le 3$, that is, $P_2$, $P_3$, and $K_3$, are clear. 

Let $n(G)\ge 4$ and assume first that $G$ is connected. Then, since $G$ is $P_4$-free, that is, a cograph, $G$ is the join of two smaller cographs $G_1$ and $G_2$. In the first subcase assume that one of $G_1$ and $G_2$ is connected, say $G_1$. Then $G_1$ is a smaller, connected cograph that satisfies the induction assumptions, hence $G_1$ is $\gamma_{tg}$-perfect with $\gamma_t(G_1) = 2$. Let Dominator start the total domination game played on $G$ with an optimal start vertex from the game played on $G_1$. If Staller replies with a move in $G_2$, the game is over. And if Staller replies with a move in $G_1$, the game is also over because $\gamma_{tg}(G_1) = 2$. In any case, $\gamma_{tg}(G) = \gamma_t(G) = 2$. In the second subcase none of $G_1$ and $G_2$ is connected. Since $G$ contains no induced $\overline{2P_3}$, we infer that at least one of $G_1$ and $G_2$, say $G_1$, is edge-less. Then the first move of Dominator on a vertex from $G_1$ forces Staller to play her first move on $G_2$, so that $\gamma_{tg}(G) = \gamma_t(G) = 2$ holds again.

If $G$ is not connected, then it consists of $k\ge 2$ components, each being an isolate- and  $(P_4, \overline{2P_3})$-free graph. Then by the above, for each component $G'$ we have $\gamma_{tg}(G') = \gamma_t(G') = 2$. Taking into account that each connected component is a co-graph and hence a join, it is now straightforward that $\gamma_{tg}(G') = \gamma_t(G') = 2k$.   

We have thus proved that if $G$ is an isolate- and $(P_4, \overline{2P_3})$-free graph, then $\gamma_t(G) = \gamma_{tg}(G)$. Since each induced subgraph of $G$ is also a $(P_4, \overline{2P_3})$-free graph, $G$ is $\gamma_{tg}$-perfect. 
\qed 

With similar but simpler arguments as used in the proof of Proposition~\ref{prp:gamma-tg-perfect} we also get the following. 
 
\begin{proposition}
 	An isolate-free  graph $G$ is $\gamma_{tg}'$-perfect, if and only if it is a cograph.
 \end{proposition}

%%%%%%%%%%%%%%%%%%%%%%%%%%%%%%%%%%%%%%%%%%%%%%%%%%%
%%%%%%%%%%%%%%%%%%%%%%%%%%%%%%%%%%%%%%%%%%%%%%%%%%%
\section{Concluding remarks}
\label{sec:remarks}
%%%%%%%%%%%%%%%%%%%%%%%%%%%%%%%%%%%%%%%%%%%%%%%%%%%
%%%%%%%%%%%%%%%%%%%%%%%%%%%%%%%%%%%%%%%%%%%%%%%%%%%

By a computer search we have obtained all $\gamma_g$-perfect and all minimally $\gamma_g$-imperfect graphs on up to $9$ vertices. The results are presented in Table~\ref{tbl:perfect}.
 
 \begin{table}[h!!]
 	\begin{center}
 		\begin{tabular}{c || c | c | c}
 			$n$ & perfect & perfect &   min.\ imperfect\\
 			 & all & connected &\\ \hline \hline
 			3 & 4 & 2 & 0\\
 			4 & 11 & 6 & 0\\
 			5 & 32 & 19 & 2\\
 			6 & 122 & 81 & 8\\
 			7 & 536 & 386 & 1\\
 			8 & 2754 & 2102 & 1\\
 			9 & 15752 & 12476 & 1\\
 		\end{tabular}
 	\end{center}
 	\caption{The number of all, and the number of connected $\gamma_g$-perfect graphs on $n$ vertices. The last column represents the number of all minimally $\gamma_g$-imperfect graphs on $n$ vertices.}
 	\label{tbl:perfect}
 \end{table}

In Proposition~\ref{prop:min-imperfect} we have proved that the anti-holes $\overline{C_n}$, $n\ge 5$, are minimally $\gamma_g$-imperfect graphs and also identified eight additional examples of  minimally $\gamma_g$-imperfect graphs.   Using computer we have verified that there are no additional such graphs on up to nine vertices. Based on these facts, as well as on the results of Theorems~\ref{thm:connected-imperfect} and~\ref{thm:triangle-free-imperfect}, we pose:

\begin{conjecture}
	There are no other minimally $\gamma_g$-imperfect graphs but those listed in Proposition~\ref{prop:min-imperfect}.  
\end{conjecture}

\section*{Acknowledgements}

We acknowledge the financial support from the Slovenian Research Agency (research core funding No.\ P1-0297 and projects J1-9109, j1-1693, N1-0095, N1-0108).

\end{document}